\documentclass[11pt]{article}

\usepackage[T1]{fontenc}
\usepackage[a4paper,margin=25mm]{geometry}

\usepackage{amsmath,amssymb,amsthm}
\usepackage{mathtools}  
\usepackage{graphicx}
\graphicspath{{figures/}}  
\usepackage[labelsep=period,font=small,labelfont=bf]{caption}
\usepackage{tikz}
\usetikzlibrary{arrows.meta}
\usepackage[section]{placeins}  
\usepackage{enumitem}
\setlist[enumerate]{label=$(\arabic*)$}
\usepackage[colorlinks,linkcolor=blue,citecolor=red,urlcolor=blue]{hyperref}  
\hypersetup{
  pdftitle={Interactive torus parametrization of the one-dimensional Fibonacci quasicrystal},
  pdfauthor={Mykhailo O. Koreshkov, Maryna O. Nesterenko},
}

\newtheorem{proposition}{Proposition}

\DeclareMathOperator{\dens}{dens}
\DeclareMathOperator{\vol}{vol}
\newcommand{\R}{\mathbb{R}}
\newcommand{\Z}{\mathbb{Z}}
\newcommand{\N}{\mathbb{N}}

\newcommand{\Ep}{E_{\mathrm{phys}}}
\newcommand{\Ei}{E_{\mathrm{int}}}
\newcommand{\piphys}{\pi_{\mathrm{phys}}}
\newcommand{\piint}{\pi_{\mathrm{int}}}
\newcommand{\hdlattice}{\widetilde L}             
\newcommand{\diffmeas}[1]{\widehat{\gamma}_{#1}}  

\newcommand{\eprint}[1]{\href{https://arxiv.org/abs/#1}{arXiv:#1}}

\title{Interactive Torus Parametrization of the One-Dimensional Fibonacci Quasicrystal}

\author{
  Mykhailo O.~Koreshkov\\
  \small Institute of Mathematics of NAS of Ukraine, Kyiv, Ukraine\\
  \small \href{mailto:mykhailo.koreshkov@gmail.com}{mykhailo.koreshkov@gmail.com},
  ORCID: \href{https://orcid.org/0009-0009-1968-0580}{0009-0009-1968-0580}
  \and
  Maryna O.~Nesterenko\\
  \small University ``Kyiv School of Economics'', Kyiv, Ukraine\\
  \small Institute of Mathematics of NAS of Ukraine, Kyiv, Ukraine\\
  \small \href{mailto:maryna.nesterenko@gmail.com}{maryna.nesterenko@gmail.com},
  ORCID: \href{https://orcid.org/0000-0001-6352-2365}{0000-0001-6352-2365}
}

\date{}  

\begin{document}

\maketitle

\begin{abstract}
We present an interactive computational framework for exploring the torus parametrization of the one-dimensional Fibonacci quasicrystal and its local hull. 
While cut-and-project schemes and Fourier analysis of local functions have been actively explored, a purely algebraic treatment still requires deeper geometric intuition regarding the underlying dynamical systems to advance these methods further.
We develop an open-source interactive Python/Jupyter tool that links coordinates in the fundamental domain of the embedding lattice directly to the physical tiling and word generation, offering a transparent visual analysis of translation orbits and singular boundary cases.
The tool is the basis for our further computational work, and it can also be used in teaching and to illustrate the geometry of the local hull.

\medskip\noindent
\textit{Keywords:} mathematical quasicrystals; cut-and-project; local hull; torus parametrization; Python.

\noindent
\textit{2020 Mathematics Subject Classification:} 52C23 (primary); 37B52, 68U05 (secondary).
\end{abstract}

\section{Introduction}\label{sec:intro}

The theory of mathematical and physical quasicrystals is a dynamic and multifaceted field linking harmonic analysis, discrete geometry, and aperiodic dynamical systems. While the structural foundations of regular model sets and pure point diffraction are firmly established, computational tools for the interactive exploration of their local hulls and function spaces remain in high demand. 
From the mathematical point of view, quasicrystals are point sets without translational symmetry but with a certain self-similarity, especially at large scales, which manifests itself in their geometric and spectral properties.
Mathematical quasicrystals are also related to dynamical systems and to tilings of space with finitely many local configurations~\cite{baake2002guide, baake2007characterization}.

In this paper, we focus on the torus parametrization of the local hull of the one-dimensional Fibonacci quasicrystal. Building upon the theoretical and computational framework of almost periodic functions introduced by Moody et al.~\cite{moody2008computing}, we present an interactive Python implementation designed to visualize the dynamics on the torus $\mathbb{T}^2$, track the generation of tile sequences, and inspect the structural bifurcations occurring at singular points.

An important property of crystalline materials is their ability to scatter radiation and produce discrete diffraction patterns that reflect the internal structure of the sample.
Traditionally, such diffraction was assumed to be exclusive to periodic structures, such as lattices and, in particular, Bravais lattices.
Crystal lattices are subject to a strict limitation on the admissible rotational symmetries, imposed by the crystallographic restriction theorem.

The discovery of physical quasicrystals (materials with ``forbidden'' symmetries and discrete diffraction patterns)
strongly stimulated the development of the new field and created the need for mathematical models of quasicrystalline structures.
The first studies in this field were devoted to almost periodic functions, i.e., trigonometric polynomials with irrational ratios of periods.
Such a function can be viewed as the restriction of a function of a higher-dimensional argument to a lower-dimensional subspace, which corresponds to a projection at an irrational angle~\cite{Bohr}.

The definition of a mathematical quasicrystal depends on the context, but it is generally understood to be a point set with a discrete diffraction pattern.
Such point sets model the positions of atoms in physical quasicrystals.

A modern approach to constructing quasicrystals is the cut-and-project method, in which a mathematical quasicrystal is obtained as a projection of a higher-dimensional lattice onto a lower-dimensional subspace.
The initial ideas of this approach were formulated by Meyer in harmonic analysis~\cite{meyer1972algebraic}, and later Moody and Patera~\cite{moody-patera} adapted and developed the method for quasicrystals.

To construct a quasicrystal, one can take the lattice $\hdlattice = \Z^{d+p} \subset \R^d \times \R^p$,
choose a physical subspace $\Ep\cong \R^d$ at an irrational angle to the lattice,
and an internal subspace $\Ei\cong \R^p$ orthogonal to the physical one.
Next, we choose a window $W \subset \Ei$ and project the elements $\tilde x\in \hdlattice$
onto the physical space, $x = \piphys(\tilde x)$, and onto the internal space, $x' = \piint(\tilde x)$.
In a more general version of the cut-and-project construction, a root lattice can be used instead of $\hdlattice = \Z^{d+p}$.
A typical construction is summarized in the following diagram, where hooked arrows denote inclusions:
\begin{equation*}
\begin{tikzpicture}[
  baseline=(current bounding box.center),
  hookarrow/.style={arrows={Hooks[right, scale=1.5, sep=4pt]-Stealth[sep=4pt]}},
  hookarrowl/.style={arrows={Hooks[left, scale=1.5, sep=4pt]-Stealth[sep=4pt]}},
  map/.style={arrows={_[sep=4pt].-Stealth[sep=4pt]}}
]
\node (R1) at (0,0) {$\R^d$};
\node (R2) at (3,0) {$\R^d\times\R^p$};
\node (R3) at (6,0) {$\R^p$};
\node (W)  at (8.2,0) {$W$};
\node (L0) at (-2.5,-1.75) {$\Lambda(W)$};
\node (L1) at (0,-1.75) {$L$};
\node (L2) at (3,-1.75) {$\hdlattice$};
\node (L3) at (6,-1.75) {$L'$};
\draw[map] (R2) -- node[below] {$\piphys$} (R1);
\draw[map] (R2) -- node[below] {$\piint$} (R3);
\draw[hookarrowl] (W) -- (R3);
\draw[hookarrow] (L0) -- (L1);
\draw[hookarrow] (L1) -- (R1);
\draw[hookarrow] (L2) -- (R2);
\draw[hookarrow] (L3) -- (R3);
\draw[map] (L2) -- node[above, font=\scriptsize] {1-1} (L1);
\draw[map] (L2) -- node[above, font=\scriptsize] {dense image} (L3);
\end{tikzpicture}
\end{equation*}

A \emph{model set} (quasicrystal) $\Lambda$ is defined as
\begin{equation*}
  \Lambda = \Lambda(W) = \{x \in L \mid x' \in W\} = \piphys \bigl(\hdlattice \cap \piint^{-1}(W)\bigr).
\end{equation*}
We consider only \emph{regular} model sets~\cite{baake2013aperiodic}, i.e., those for which
$\piint(\hdlattice)$ is dense in $\Ei$,
$\piphys|_{\hdlattice}$ is injective,
and the window $W$ is a bounded Borel set with nonempty interior and boundary of measure zero.
In this paper, we only consider windows that are boxes, $W = \prod_{i=1}^p[a_i,b_i]$, although in general a window can have a complicated, fractal structure~\cite[Fig.~5.2]{baake2002guide}.

Regular model sets are uniformly discrete and have a discrete diffraction pattern~\cite{baake2013aperiodic}.
Let us examine these notions in more detail.

Let $X \subset \R^d$ be a uniformly discrete point set, let $\delta_x$ be the Dirac delta (distribution) at $x$,
and let $\omega_X = \sum_{x \in X} \delta_x$ be the \emph{Dirac comb}, the measure associated with $X$.
Uniform discreteness ($\inf_{x \neq y \in X} \|x-y\| > 0$) guarantees that $\omega_X$ is a locally finite measure.

Let $B_r(x) = \{y \mid \|y-x\| < r\}$.
We define the \emph{autocorrelation measure} of $X$ as
\[
  \gamma_X = \lim_{r\to\infty} \frac{1}{\vol(B_r(0))} \sum_{x,y \in X \cap B_r(0)} \delta_{x-y}.
\]
The \emph{diffraction measure} of $X$ is the measure $\diffmeas{X}$, the Fourier transform of the autocorrelation,
defined by its action on Schwartz functions $f \in \mathcal{S}(\R^d)$ as $\diffmeas{X}(f) = \gamma_X\bigl(\widehat f \bigr)$.
If $X$ is a lattice or a model set, then both measures exist, and the autocorrelation has the form $\gamma_X = \sum_{z \in \Delta_X}\nu(z)\delta_z$,
where $\nu(z)$ are real coefficients and $\Delta_X = \{x-y \mid x,y \in X\}$~\cite{baake2002guide}.

Consider a lattice $\Gamma = \Z b_1 \oplus \dots \oplus \Z b_d$ with linearly independent basis vectors $b_1, b_2, \dots, b_d$.
Let $\Gamma^* = \{y \in \R^d \mid \langle y, x\rangle \in \Z \ \forall x \in \Gamma\}$ be the dual lattice, and let
$\dens(\Gamma) = 1/\lvert\det(b_1, b_2,\dots,b_d)\rvert$.
Then $\Delta_\Gamma = \Gamma$, the autocorrelation measure $\gamma_\Gamma$ exists,
and the diffraction measure is discrete: $\diffmeas{\Gamma} = \dens(\Gamma)^2 \sum_{y\in\Gamma^*}\delta_y$.

Point sets with similar properties are \emph{Delone sets}: uniformly discrete sets that are also \emph{relatively dense},
i.e., there exist $0< R,r < \infty$ such that $\forall x \in \R^d\ \exists\, y \in X \colon \|x-y\| < R$ and
$\forall x,y \in X,\ x \neq y \colon \|x-y\| > r$.
In other words, Delone sets are discrete sets in which the distance from any point to its nearest neighbor is
bounded from above and from below, and there are no gaps larger than $R$.
In the language of metric geometry, a Delone set is an $r$-separated $R$-net.

We call Delone sets with a discrete diffraction pattern, i.e., with a pure point diffraction measure, \emph{mathematical quasicrystals}.
This is expressed by the following conditions:
\begin{gather*}
  \exists\ \dens(X) < \infty, \quad \dens(X) \coloneqq \lim_{R \to \infty} \frac{\#(X \cap B_R(0))}{\vol(B_R(0))}, \\
  \exists\ \diffmeas{X}, \quad \diffmeas{X} \text{ is a pure point measure}.
\end{gather*}

Importantly, the diffraction pattern of a quasicrystal is discrete but \emph{dense},
and it has a self-similar structure if one considers only the diffraction peaks whose intensity exceeds a given threshold.
There are also other ways to construct quasicrystals, from the point of view of dynamical systems,
for example, by substitution rules, subdivision rules, or finite automata~\cite{tatham2025finite}.

Higher-dimensional quasicrystals and the corresponding tilings are obtained by generalizing the construction
of the one-dimensional Fibonacci quasicrystal (see Section~\ref{sec:fibonacci}):
the hypercubes of the lattice that fall within the window are projected onto the physical subspace.
Many algorithms exist for constructing point sets and tilings in a two-dimensional physical space.
Popular algorithms use the traces of the intersections of lattice hyperplanes with the physical plane to find admissible points,
for example, the generalized dual method~\cite{socolar1985gendualmethod},
de Bruijn's method~\cite{debruijn, gregegantechnical},
and the multigrid method~\cite{gglouser}.

Conversely, one may need to reconstruct a quasicrystal from its diffraction pattern obtained from experimental data.
A well-known tool for this is Superflip~\cite{palatinus2007superflip},
which solves crystal structures in arbitrary dimension with the iterative charge-flipping algorithm.
Another approach is to process surface microscopy images with machine learning and computer vision algorithms,
see, e.g.,~\cite{kender2025automatic} and~\cite{chang2023}.
Given noisy coordinates of the points of a model set or a higher-dimensional lattice, one may need to determine its symmetries;
methods for this are formulated in the classic paper~\cite{lepage1987computer} and in the more recent~\cite{zwart2006exploring}.
More recently, the papers~\cite{petrache2023almost,ouguz2017hyperuniformity,yakir2022recovering} established conditions
under which a randomly perturbed structure can be recovered almost surely.

Note that model sets and tilings can be built not only on hypercubic lattices but also on root lattices, i.e., lattices generated by root systems.
For some quasicrystals, this gives a simpler description and a lower dimension of the ambient space.
At the same time, the symmetry operators become more complicated, and their description requires the theory of Coxeter and Weyl groups.
For example, the Tübingen triangle tiling can be obtained from the root lattice~$A_4$ by projecting its Delaunay cells,
and the Penrose tiling by projecting its Voronoi cells~\cite{baake1990planar}.

Of central importance are local functions, whose values at any point are determined by the local patch of the quasicrystal around that point. In~\cite{moody2008computing} it was shown that every bounded local function $f$ on a regular model set lifts to a function $\widetilde{F}$ on the torus $\mathbb{T}^{2d}$ via the parametrization of the local hull, which can then be approximated uniformly by finite Fourier series. Our software directly demonstrates the geometric support of such functions on the fundamental domain, providing visual intuition for how discontinuous functions on the torus generate almost periodic signals along irrational winding lines.


\section{The Fibonacci quasicrystal}\label{sec:fibonacci}

The classic and arguably the most studied example of a quasicrystal is the family of one-dimensional \emph{Fibonacci quasicrystals}.
The specific form of the model set depends on the choice of the lattice, the physical and internal subspaces, and the window.

One construction takes the square lattice $\hdlattice=\Z^2$ and places the physical subspace
at an irrational angle to it, for example, $\arctan \tau^{-1}$, where $\tau=(1+\sqrt{5})/2$ is the golden ratio.
Let $e_{\mathrm{phys}} = (1, \tau^{-1})$ and $e_{\mathrm{int}} = (-\tau^{-1}, 1)$ be orthogonal vectors
along the physical and internal subspaces, respectively.
Instead of $\tau$, one can take the positive root of $x^2=m x+1$ for another value of~$m\in\N$.

We choose the window $W$ as the projection onto the internal subspace of the Voronoi cell of the origin of $\Z^2$,
i.e., $W = \piint(\{x \in \R^2 \mid \|x\|_\infty < 1/2\})$.
Note that this window is symmetric about the origin.

\begin{figure}[htbp]
  \centering
  \includegraphics[
    width=0.5\linewidth,
    alt={A square lattice of black dots with a line of slope about 0.62 through the origin,
         a light blue strip around this line, and a perpendicular line.
         The lattice points inside the strip form a staircase path,
         and their projections onto the line are shown as red dots.
         A shaded unit square at the origin marks the fundamental cell.
         The path takes horizontal and vertical steps in a seemingly random order.
         The gaps between the projected red dots have two distinct lengths:
         a horizontal step projects to a long gap, and a vertical step to a short gap.}
  ]{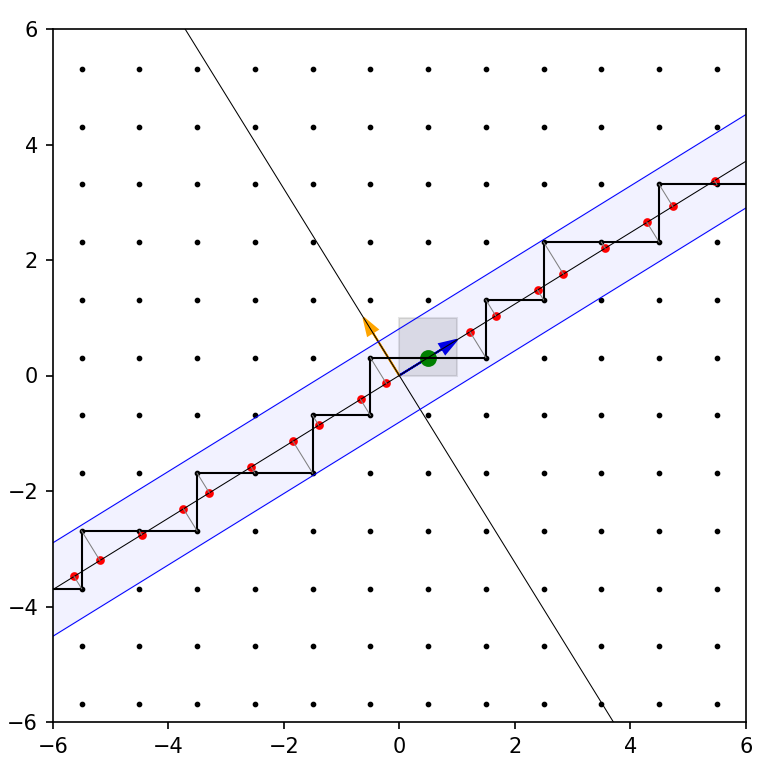}
  \caption{The Fibonacci quasicrystal as the projection of the lattice $\Z^2$
    onto the physical line $y=\frac{2}{1+\sqrt{5}} x$ with slope $\tau^{-1} = (\sqrt{5}-1)/2 \approx 0.618$ $(\frac{\sqrt{5}-1}{2} = \tau - 1 = \tau^{-1})$.}
  \label{fig:classic_fibonacci_1}
\end{figure}

As a result, we obtain a nonperiodic discrete set $\Lambda$ (see Fig.~\ref{fig:classic_fibonacci_1}).
Neighboring points of $\Lambda$ can be at only two different distances from each other, which defines the \emph{tile type}.
Thus, we only consider quasicrystals that are tilings of $\R$ by two types of tiles, ``long'' and ``short'',
which correspond to a step to the right and a step up on the lattice, respectively (singular cases aside).

\section{The local hull}\label{sec:hull}

The definition of a Delone set given in Section~\ref{sec:intro} can be restated equivalently:
a set $\Lambda \subset \R^d$ is a Delone set if there exist $r, R > 0$ such that
\[
  \inf\bigl\{\|x - y\| \mid x,y \in \Lambda,\ x \neq y\bigr\} \ge r \quad \text{and} \quad \sup_{x \in \R^d}\,\inf_{y \in \Lambda}\|x - y\| \le R.
\]
Let $t \in \R^d$ act on $\Lambda$ by translation (shift), i.e., $t + \Lambda = \{t + x \mid x \in \Lambda\}$.

Delone sets $A$ and $B$ are called \emph{locally indistinguishable} (LI) if every \emph{patch} of $A$, i.e.,
the finite configuration $A \cap B_r(x)$ of the points of $A$ within distance $r$ of a point $x \in A$,
occurs in $B$ up to translation, and vice versa~\cite{baake2002guide}.
Formally, $\forall x\in A\ \forall r>0\ \exists t \in \R^d \colon t + (A\cap B_r(x)) = B \cap B_r(t + x)$,
and the same holds with $A$ and $B$ interchanged.
Thus, local indistinguishability is an equivalence relation.
Moreover, Delone sets are equipped with the \emph{local topology}~\cite{moody2008computing},
in which two sets are close if they agree, up to a small error $\varepsilon$, within a large radius $R$;
such sets are called $(R,\varepsilon)$-close.
The local topology is defined by the base of open sets of the form
\begin{equation*}
  U_{R,\varepsilon}(A) = \left\{ B \  \middle|
    \begin{array}{l}
      B \cap B_R(0) \subset A + B_\varepsilon(0), \\
      A \cap B_R(0) \subset B + B_\varepsilon(0)
    \end{array}
    \right\}.
\end{equation*}

The \emph{local hull} of a Delone set $\Lambda$ is the topological space $X(\Lambda)$ obtained
as the closure of the set of translates of $\Lambda$ in the local topology:
\[
  X(\Lambda) = \overline{\{t + \Lambda \mid t \in \R^d\}}.
\]
A tiling, or a Delone set, has \emph{finite local complexity} if, for every $r > 0$,
it has only finitely many different patches of radius $r$ up to translation.

\begin{proposition}[{\cite[Lemma~1]{radin1992}}]\label{prop:lemma1}
  The local hull of a tiling of the plane with finite local complexity (e.g., a Penrose tiling) consists
  of all tilings built from the same set of local configurations.
\end{proposition}

\begin{proposition}[{\cite[Lemma~2]{radin1992}}]
  The topological space $X(\Lambda)$ is compact, and the translation action of $\R^d$ on $\Lambda$ lifts
  to a continuous translation action on $X(\Lambda)$.
\end{proposition}

By Proposition~\ref{prop:lemma1}, the local hull consists of all Delone sets that are locally indistinguishable from $\Lambda$.
For this reason, the local hull is also called the \emph{local indistinguishability class} (LI class).
For the lattice $\Gamma = \Z^d \subset \R^d$, the local hull $X(\Gamma)$ is the torus $\mathbb{T}^d = \R^d / \Z^d$.
All sets locally indistinguishable from $\Gamma$ are translates of $\Gamma$,
and they are parametrized by the points of the fundamental cell with opposite faces identified.

The realizations of the Fibonacci quasicrystal with a given window and slope are parametrized in a similar way,
by the torus $\mathbb{T}^2$ on the fundamental cell.
Let us revisit the cut-and-project scheme, now with physical and internal subspaces of equal dimension, $d=p$:
\[
  \Lambda = \Lambda(W) = \{x \in L \mid x' \in W\} = \piphys\bigl(\hdlattice \cap \piint^{-1}(W)\bigr),
\]
where $\hdlattice \subset \R^d\times \R^d$ is a lattice, $L\subset \Ep \cong \R^d$, and $L' \subset \Ei \cong \R^d$.

For $(a,b) \in \R^d \times \R^d$, denote
\[
  \Lambda_{(a,b)} = a + \Lambda(-b + W) = \{x+a \mid x \in L,\ x'+b \in W\}.
\]
Vectors $a$ and $b$ from $\R^d$ act as translations in the physical and internal spaces, respectively.
If $(a,b) - (a',b') \in \hdlattice$, the lattice in $\R^d\times\R^d$, then $\Lambda_{(a,b)} = \Lambda_{(a',b')}$, so $\Lambda_{(a,b)}$ depends only on the equivalence class $(a,b)_L \coloneqq (a,b) + \hdlattice$, a point of the torus $(\R^d\times\R^d)/\hdlattice \cong \mathbb T^{2d}$.
For $\hdlattice = \Z^2$, writing $a + b \in \R^2$ in the standard basis and reducing both coordinates modulo~1 identifies
$(a,b)_L$ with a point of the fundamental cell $[0,1)^2$.

All further analysis of functions associated with quasicrystals, as well as the construction of nonperiodic discrete sets
with a large number of points, relies on Propositions~\ref{prop_moody} and~\ref{prop_s_b},
proved in~\cite{moody2008computing} and~\cite{schlottmann2000generalized, baake1997torus}, respectively.

\begin{proposition}\label{prop_moody}
  The set $\{(t,0)_L \mid t\in\R^d\}$ is dense in the torus $(\R^d\times\R^d)/\hdlattice$.
\end{proposition}
Intuitively, this means that the points $(t,0)_L$ trace a line of irrational slope on the torus, and this line fills the torus densely.

\begin{proposition}\label{prop_s_b}
  There exists a \emph{torus parametrization of the local hull}, i.e., a continuous map
  \[\beta \colon X(\Lambda) \longrightarrow \mathbb T^{2d}\]
  such that
  \begin{enumerate}
    \item $\beta$ is surjective;
    \item $\beta$ is injective almost everywhere (with respect to the Haar measure on the torus);
    \item $\beta(t+\Lambda') = \beta(\Lambda') + (t, 0)_L$ for all $t \in \R^d$ and all $\Lambda' \in X(\Lambda)$;
    \item $\beta(t + \Lambda) = (t, 0)_L$ for all $t \in \R^d$.
  \end{enumerate}
\end{proposition}

A point of the torus is called \emph{singular} if $\beta$ is not injective over it, i.e., if several elements of $X(\Lambda)$ are mapped to this point.
For example, this happens when a point of the translated lattice falls exactly on the boundary of the strip $\piint^{-1}(W)$: this point can either be included or excluded, which produces two different Fibonacci chains (see Fig.~\ref{fig:fibonacci-interface}).
In other words a point of the torus is singular whenever a point of the shifted lattice lies precisely on the boundary $\partial \pi_{\mathrm{int}}^{-1}(W)$. 
Such points represent boundary ambiguities where the torus map is not injective: two distinct quasicrystal configurations map to the same point on $\mathbb{T}^2$. 
In harmonic analysis, this manifests as a localized Gibbs-like misfit or `fuzziness' of Fourier approximations near the origin. Our interactive notebook explicitly allows users to test these boundary values (e.g., $x=0.5, y=0.5$) and observe the non-uniqueness of the resulting tiling.

\section{Interactive software}\label{sec:software}

To illustrate the torus parametrization of the local hull for the Fibonacci quasicrystal,
we developed interactive Python software in the form of a Jupyter notebook;
its source code is freely available on GitHub.\footnote{\url{https://github.com/mkrooted256/fibonacci-torus-notebook/}}

In each figure, the semi-transparent strip is the strip around the physical space that projects onto the window $W$ in the internal space.
The shaded square is the fundamental cell of the lattice, which parametrizes the local hull as a torus.

The parametrization by $(x,y)_L$ is realized by translating the lattice by the vector $(x,y)$,
while the window and the subspaces stay fixed.
This corresponds to the parametrization by $(a,b)_L$ above, with $(x,y) = a + b$, $a \in \Ep$, $b \in \Ei$.

The software implements the parametrization in $xy$ and $uv$ coordinates:
$xy$ are coordinates in the plane restricted to the fundamental domain $[0,1]^2$,
while $uv$ are coordinates with respect to the basis vectors of the physical and internal spaces.

To facilitate reproducibility and direct comparison with the computational implementation in our repository, we explicitly specify the coordinate transformation between the standard torus coordinates $(x, y) \in [0, 1)^2$ and the physical/internal decomposition $(u, v) \in E_{\mathrm{phys}} \times E_{\mathrm{int}}$. Since the basis vectors $e_{\mathrm{phys}} = (1, \tau^{-1})^{\top}$ and $e_{\mathrm{int}} = (-\tau^{-1}, 1)^{\top}$ are mutually orthogonal with equal norms $\|e_{\mathrm{phys}}\|^2 = \|e_{\mathrm{int}}\|^2 = 1 + \tau^{-2} = 2 + \tau^{-1}$, any point can be decomposed as
\begin{equation}\label{eq:coord_transform}
\begin{pmatrix} x \\ y \end{pmatrix}
= u \frac{e_{\mathrm{phys}}}{\|e_{\mathrm{phys}}\|^2} + v \frac{e_{\mathrm{int}}}{\|e_{\mathrm{int}}\|^2}
= \frac{1}{1 + \tau^{-2}}
\begin{pmatrix}
1 & -\tau^{-1} \\
\tau^{-1} & 1
\end{pmatrix}
\begin{pmatrix} u \\ v \end{pmatrix}.
\end{equation}
Conversely, the projection coordinates $(u, v)$ evaluated in the algorithms are obtained via the direct orthogonal projections $u = \langle (x, y)^{\top}, e_{\mathrm{phys}} \rangle = x + \tau^{-1} y$ and $v = \langle (x, y)^{\top}, e_{\mathrm{int}} \rangle = -\tau^{-1} x + y$. This explicit matrix representation aligns the geometric projection steps directly with the vectorised transformations used throughout the accompanying Jupyter notebook.

\begin{figure}[htbp]
  \centering
  \includegraphics[
    width=0.38\linewidth,
    alt={Cut-and-project picture for u = 0.2, v = 0:
         lattice, strip around the physical line, staircase of lattice points inside the strip,
         their projections as red dots, and the fundamental cell with the parameter point marked in green.}
  ]{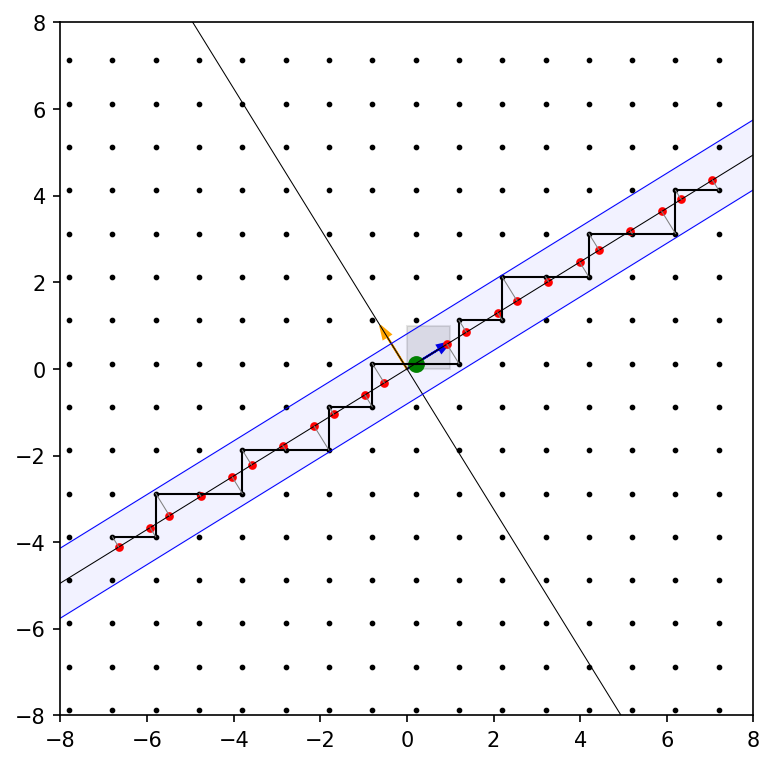}%
  \hfill
  \includegraphics[
    width=0.38\linewidth,
    alt={The same picture for u = 0.5, v = 0:
         the red points are shifted along the physical line, and their pattern is unchanged.}
  ]{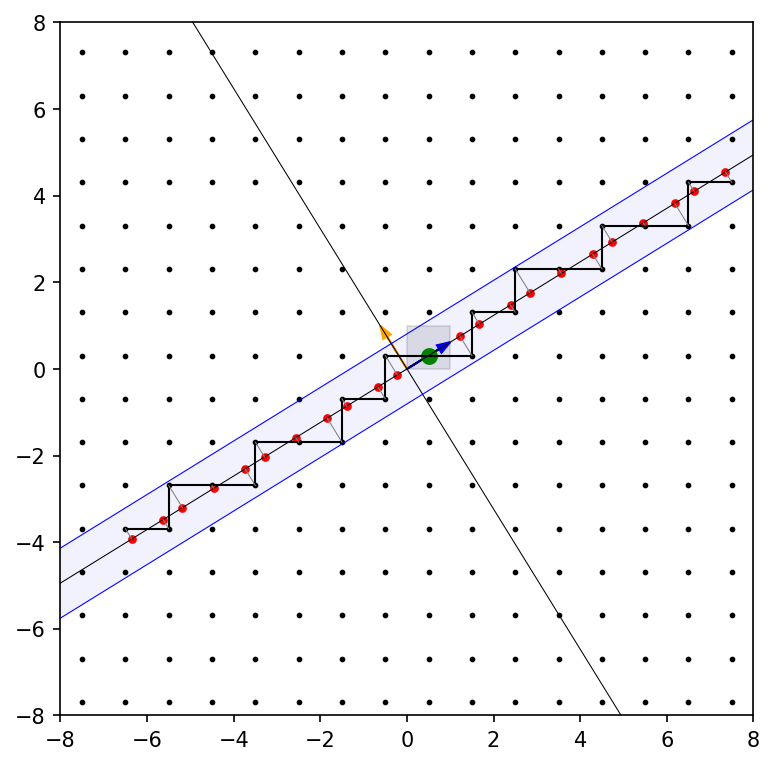}
  \caption{Left: $u = 0.2$, $v = 0$; right: $u = 0.5$, $v = 0$.
    The configuration is the same, translated in physical space.}
  \label{fig:fibonacci-phys}
\end{figure}

\begin{figure}[htbp]
  \centering
  \includegraphics[
    width=0.38\linewidth,
    alt={Cut-and-project picture for u = 0, v = 0,
         with the parameter point at the corner of the fundamental cell.}
  ]{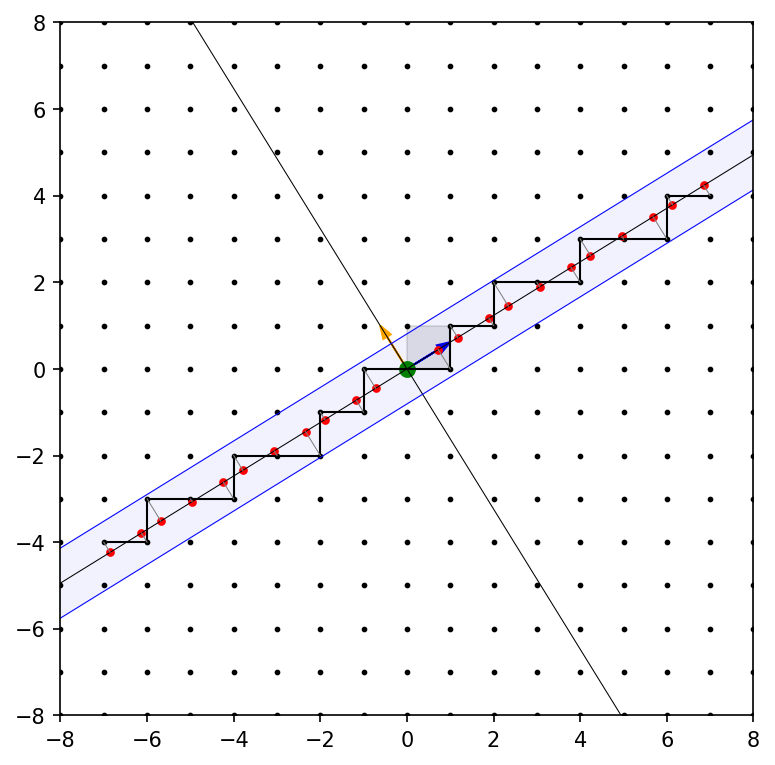}%
  \hfill
  \includegraphics[
    width=0.38\linewidth,
    alt={The same picture for u = 0, v = 0.5:
         the strip selects different lattice points,
         so the sequence of red points along the physical line is different.}
  ]{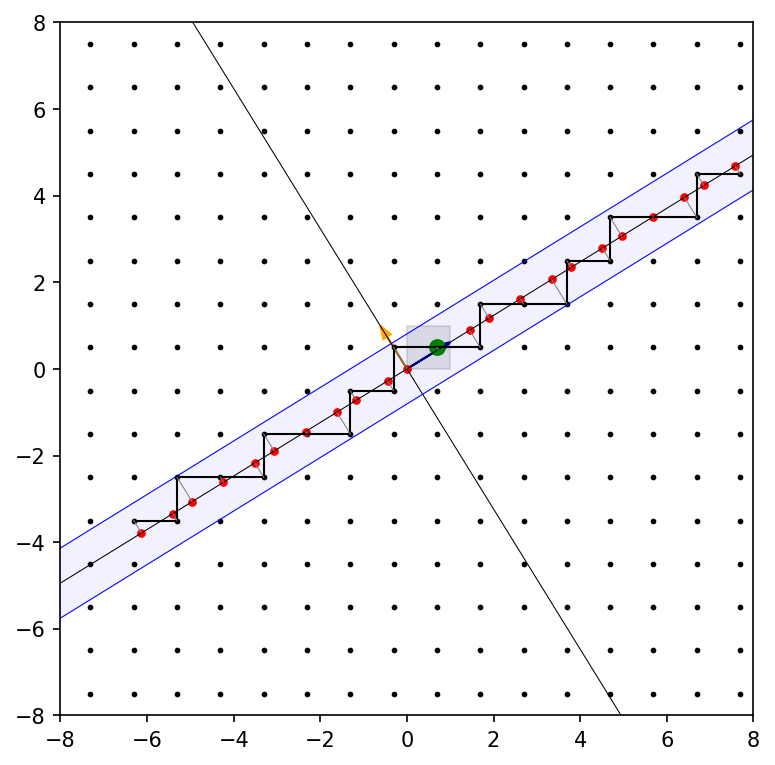}
  \caption{Left: $u = 0$, $v = 0$; right: $u = 0$, $v = 0.5$.
    The configuration changes (translation in internal space).
    The winding around the torus is visible:
    the trajectory starts in the lower left corner and continues from the lower right corner.}
  \label{fig:fibonacci-int}
\end{figure}

Changing only the $u$ coordinate gives the realizations $\Lambda_{(t,0)}$, i.e.,
it simply translates the model set in physical space (Fig.~\ref{fig:fibonacci-phys}).
Changing the $v$ coordinate changes the realization of the quasicrystal.
This is clearly visible in the sequence of tile types, which the notebook prints as an \emph{ab-word},
where $a$ stands for a long tile (a horizontal step on the lattice) and $b$ for a short tile (a vertical step):
the words differ near the origin (Fig.~\ref{fig:fibonacci-int}).
Increasing $u$ or $v$ alone makes the point that corresponds to the origin of the quasicrystal wind around the torus.

\begin{figure}[htbp]
  \centering
  \includegraphics[
    width=0.38\linewidth,
    alt={Cut-and-project picture for x = y = 0.5,
         where lattice points lie exactly on the boundary of the strip,
         so the staircase has an ambiguous zigzag near the origin.}
  ]{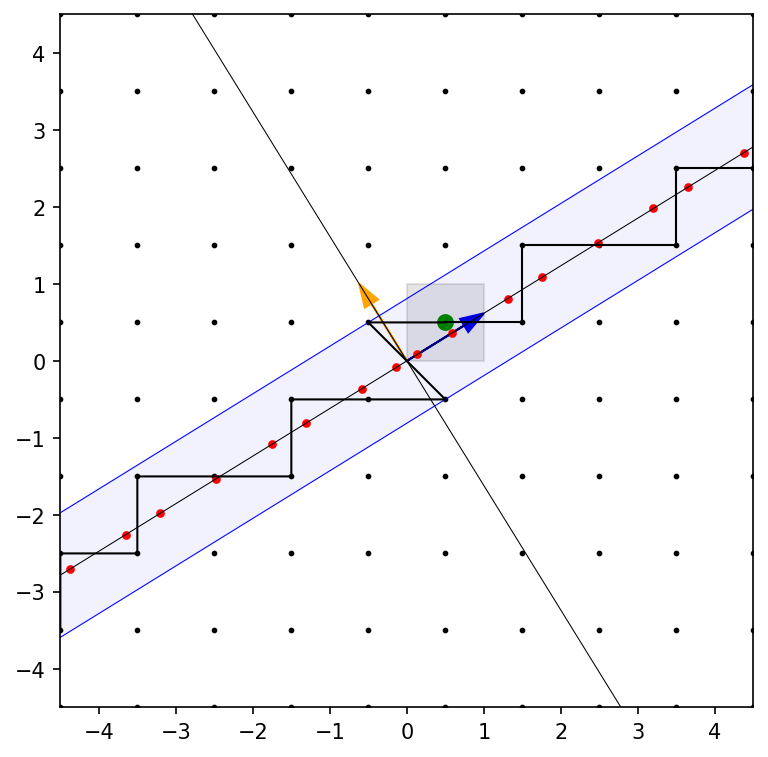}%
  \hfill
  \includegraphics[
    width=0.38\linewidth,
    alt={Screenshot of the notebook interface:
         sliders for u and v, a text readout of tau, the uv and xy coordinates, the window,
         and the tile sequence as an ab-word with a bar marking the origin,
         above the cut-and-project plot.}
  ]{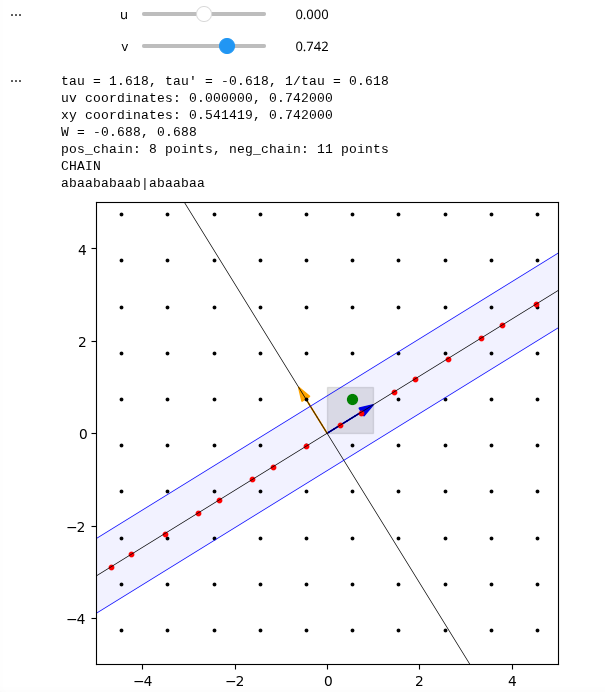}
  \caption{Left: the singular case $x = 0.5$, $y = 0.5$;
    right: the interactive interface for entering $uv$ coordinates.}
  \label{fig:fibonacci-interface}
\end{figure}

Besides the code, documentation, examples, and instructions, the Jupyter notebook provides an interactive
interface for entering $uv$ or $xy$ coordinates in the fundamental domain (Fig.~\ref{fig:fibonacci-interface}).

\section{Conclusions}\label{sec:conclusions}

We have presented the basic definitions of the theory of mathematical quasicrystals, with a focus on the Fibonacci quasicrystal and its local hull.
For a qualitative analysis and interpretation of the torus parametrization of the local hull,
we developed interactive Python software that visualizes different realizations of the Fibonacci quasicrystal and
shows how they correspond to points on the torus.

In future work, we plan to extend this interactive module to compute discrete Fourier sums of local functions in real time, implementing the refinement lattices $\widetilde{L}_N = \frac{1}{N}\mathbb{Z}[\tau]$ introduced in~\cite{moody2008computing}. 
Developing and integrating a numerical implementation of the fast Fourier transform on the torus with real-time visualization of quasicrystalline patterns will provide an integrated educational and research workbench for aperiodic structures.

\section*{Acknowledgments}

The authors thank the anonymous reviewers of the conference version~\cite{KoreshkovNesterenko2026} of this paper for their valuable comments and suggestions, which helped improve it.

\medskip
\noindent
\begin{minipage}{20mm}
  \includegraphics[
    width=20mm,
    alt={National Research Foundation of Ukraine logo}
  ]{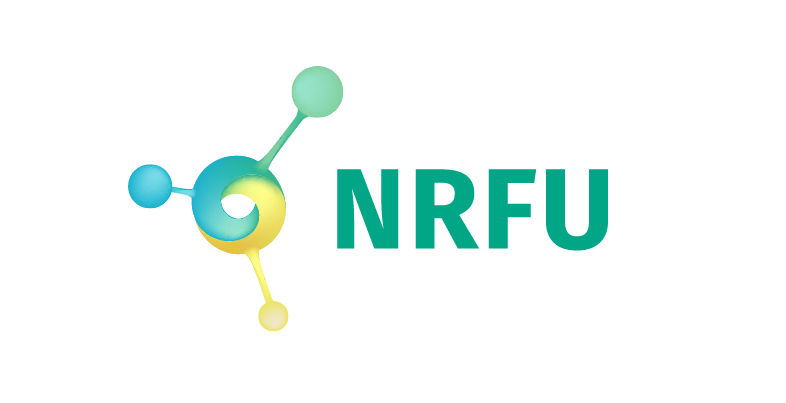}
\end{minipage}
\
\begin{minipage}{90mm}
  \footnotesize\itshape
  Prepared for publication and published with the grant support of the National Research Foundation of Ukraine within the framework of the project 2025.07/0405 ``Algebraic methods for studying equations of mathematical physics''.
  The contents of this work do not necessarily reflect the views of the National Research Foundation of Ukraine and are the sole responsibility of Institute of Mathematics of NAS of Ukraine.
\end{minipage}


\begin{thebibliography}{99}
\footnotesize\itemsep=0pt\frenchspacing

\bibitem{baake2002guide}
Baake M., \href{https://doi.org/10.1007/978-3-662-05028-6_2}{A guide to mathematical quasicrystals}, in \textit{Quasicrystals: An introduction to structure, physical properties and applications}, Springer, 2002, 17--48, \eprint{math-ph/9901014}.

\bibitem{baake2013aperiodic}
Baake M., Grimm U., \href{https://doi.org/10.1017/CBO9781139025256}{Aperiodic order. {V}ol. 1}, \textit{Encyclopedia Math. Appl.}, Vol.~149, Cambridge University Press, Cambridge, 2013.

\bibitem{baake1997torus}
Baake M., Hermisson J., Pleasants P.A.B., \href{https://doi.org/10.1088/0305-4470/30/9/016}{The torus parametrization of quasiperiodic {LI}-classes}, \textit{J.~Phys.~A} \textbf{30} (1997), no.~9, 3029--3056.

\bibitem{baake1990planar}
Baake M., Kramer P., Schlottmann M., Zeidler D., \href{https://doi.org/10.1142/S0217979290001054}{Planar patterns with fivefold symmetry as sections of periodic structures in {$4$}-space}, \textit{Internat. J. Modern Phys. B} \textbf{4} (1990), no.~15--16, 2217--2268.

\bibitem{baake2007characterization}
Baake M., Lenz D., Moody R.V., \href{https://doi.org/10.1017/S0143385706000800}{Characterization of model sets by dynamical systems}, \textit{Ergodic Theory Dynam. Systems} \textbf{27} (2007), no.~2, 341--382.

\bibitem{Bohr}
Bohr H., Almost periodic functions, Chelsea Publishing Co., New York, 1947, available at \url{https://archive.org/details/in.ernet.dli.2015.84535}.

\bibitem{debruijn}
de~Bruijn N.G., \href{https://pure.tue.nl/ws/portalfiles/portal/4344195/597566.pdf}{Algebraic theory of {P}enrose's nonperiodic tilings of the plane}, \textit{Nederl. Akad. Wetensch. Indag. Math.} \textbf{43} (1981), no.~1, 39--66.

\bibitem{gregegantechnical}
Egan G., \textit{deBruijn applet: Mathematical Details}, 2008, available at \url{https://www.gregegan.net/APPLETS/12/deBruijnNotes.html}.

\bibitem{gglouser}
Glouser G., \textit{Cut and project tiling applet on Github}, 2018, available at \url{https://gglouser.github.io/cut-and-project-tiling/}.

\bibitem{kender2025automatic}
Kender T.K., Corrias M., Franchini C., \href{https://doi.org/10.1002/aidi.202500043}{Automatic determination of quasicrystalline patterns from microscopy images}, \textit{Adv. Intell. Discov.}, 2025, \eprint{2503.05472}.

\bibitem{KoreshkovNesterenko2026}
Koreshkov M., Nesterenko~M,
\emph{Interactive torus parametrization of the one-dimensional Fibonacci quasicrystal},
Proc. Conf. Young Sci. KyivAcademUs, 2026, math.01, 9~pp.,
\href{https://doi.org/10.3842/kau.2026.math.01}{doi:10.3842/kau.2026.math.01}.

\bibitem{lepage1987computer}
Le~Page Y., \href{https://doi.org/10.1107/S0021889887086710}{{Computer derivation of the symmetry elements implied in a structure description}}, \textit{J.~Appl.~Crystallogr.} \textbf{20} (1987), no.~3, 264--269.

\bibitem{chang2023}
Liu C., Kitahara K., Ishikawa A. et al., \href{https://doi.org/10.1103/PhysRevMaterials.7.093805}{Quasicrystals predicted and discovered by machine learning}, \textit{Phys.~Rev.~Mater.} \textbf{7} (2023), 093805, 9~pages.

\bibitem{meyer1972algebraic}
Meyer Y., \textit{Algebraic numbers and harmonic analysis}, Vol.~2, Elsevier, 1972.

\bibitem{moody2008computing}
Moody R.V., Nesterenko M., Patera J., \href{https://doi.org/10.1107/S0108767308025440}{Computing with almost periodic functions}, \textit{Acta Crystallogr. Sect.~A} \textbf{64} (2008), no.~6, 654--669, \eprint{0808.1814}.

\bibitem{moody-patera}
Moody R.V., Patera J., \href{https://doi.org/10.1088/0305-4470/26/12/022}{Quasicrystals and icosians}, \textit{J.~Phys.~A} \textbf{26} (1993), no.~12, 2829--2853.

\bibitem{ouguz2017hyperuniformity}
O\u{g}uz E.C., Socolar J.E.S., Steinhardt P.J., Torquato S., \href{https://doi.org/10.1103/PhysRevB.95.054119}{Hyperuniformity of quasicrystals}, \textit{Phys.~Rev.~B} \textbf{95} (2017), 054119, 10~pages.

\bibitem{palatinus2007superflip}
Palatinus L., Chapuis G., \href{https://doi.org/10.1107/S0021889807029238}{{\textit{SUPERFLIP}---a computer program for the solution of crystal structures by charge flipping in arbitrary dimensions}}, \textit{J.~Appl.~Crystallogr.} \textbf{40} (2007), no.~4, 786--790.

\bibitem{petrache2023almost}
Petrache M., Viera R., \href{https://doi.org/10.1007/s10955-022-03059-2}{Almost sure recovery in quasi-periodic structures}, \textit{J.~Stat.~Phys.} \textbf{190} (2023), no.~2, 39, 26~pages, \eprint{2112.11613}.

\bibitem{radin1992}
Radin C., Wolff M., \href{https://doi.org/10.1007/BF02414073}{Space tilings and local isomorphism}, \textit{Geom. Dedicata} \textbf{42} (1992), no.~3, 355--360.

\bibitem{schlottmann2000generalized}
Schlottmann M., Generalized model sets and dynamical systems, in \textit{Directions in mathematical quasicrystals}, \textit{CRM Monogr. Ser.}, Vol.~13, American Mathematical Society, Providence, RI, 2000, 143--159.

\bibitem{socolar1985gendualmethod}
Socolar J.E.S., Steinhardt P.J., Levine D., \href{https://doi.org/10.1103/PhysRevB.32.5547}{Quasicrystals with arbitrary orientational symmetry}, \textit{Phys.~Rev.~B} \textbf{32} (1985), 5547--5550.

\bibitem{tatham2025finite}
Tatham S., \eprint{2512.16595}.

\bibitem{yakir2022recovering}
Yakir O., \href{https://doi.org/10.1093/imrn/rnaa316}{Recovering the lattice from its random perturbations}, \textit{Int. Math. Res. Not. IMRN} \textbf{2022} (2022), no.~8, 6243--6261, \eprint{2002.01508}.

\bibitem{zwart2006exploring}
Zwart P., Grosse-Kunstleve R., Adams P., \href{https://escholarship.org/uc/item/7ch589hv}{Exploring metric symmetry}, \textit{Lawrence Berkeley National Laboratory}.

\end{thebibliography}
\end{document}